\documentclass [12pt]{amsart}
\usepackage {amssymb}
\usepackage {amsmath}

\begin{document}

\begin{center}
\textbf{Transversality Conditions for Higher Order Infinite Horizon Discrete Time Optimization Problems}
\end{center}

\begin{center}

\end{center}
\begin{center}
Dapeng CAI $^{1}$ and Takashi Gyoshin NITTA $^{2}$
\end{center}

\begin{center}
$^{1 }$\textit{Corresponding author. Institute for Advanced Research, Nagoya University, Furo-cho, Chikusa-ku, Nagoya, 464-8601, Japan;}$^{ 2}$\textit{ Department of Mathematics, Faculty of Education, Mie University, Kurimamachiya 1577, Tsu, 514-8507, Japan}
\end{center}

\begin{center}

\end{center}
\begin{center}
\textbf{Abstract}
\end{center}

In this paper, we examine higher order difference problems:

\noindent 
$\mathop 
{\mbox{max}}\limits_{\rm {\bf c}}   
\sum\limits_{t=0}^\infty 
U( 
{{\rm {\bf c}}( t ),
{\rm {\bf c}}(t+1),\cdots ,
{\rm {\bf c}} (t+N-1),t)} $.
 Using the ``squeezing'' argument, 
we derive both Euler's condition and the transversality condition. In order 
to derive the two conditions, two needed assumptions are identified. A 
counterexample, in which the transversality condition is not satisfied 
without the two assumptions, is also presented.


\newpage
\section{Introduction}
In this paper, we consider the following reduced form model
\begin{equation}
\begin{cases}
\begin{aligned}
&\displaystyle\max_{\mathbf{c}}\displaystyle\sum_{t=0}^{\infty}
U\left(\mathbf{c}(t),\mathbf{c}(t+1),\cdots,\mathbf{c}(t+N-1),t\right)
\\
&\mathrm{subject\; to}\; 
\mathbf{c}(0)=\mathbf{c}_{0},\\
&\forall t\geq 0,\; 
\left(\mathbf{c}(t),\mathbf{c}(t+1),\cdots,\mathbf{c}(t+N-1)\right)\in X(t)\subset\left(\mathbb{R}^{n}\right)^{N},
\end{aligned}
\end{cases}
\label{eq1}
\end{equation}

\noindent where $N\in \mathbb{N}$, $U$ is a real-valued $N$th-order continuously 
differentiable function, and ${\rm {\bf c}}\equiv \left( {c_1 ,c_2 ,\cdots 
,c_n } \right)$ is $N$th-order continuously differentiable.
\footnote{ Normally, 
$U$ is defined on $\left( {\mathbb{R}^n} \right)^N\times \mathbb{R}$. 
The domain of $U$ is denoted by $X\left( t \right)$, in included in $\left( 
{\mathbb{R}^n} \right)^N$, for all $t$.}
 Notice that the objective functional of (\ref{eq1}) can be infinite. \cite{ref7} considers the 
continuous time first order differential problems: $v\left( {{\rm {\bf 
x}}\left( t \right),{\rm {\bf \dot {x}}}\left( t \right),t} \right)$. It 
generalizes the results of \cite{ref5,ref6,ref10,ref12}.
 So far, the most general form of the transversality 
conditions for continuous time version of problem (\ref{eq1}) is presented in 
\cite{ref11}, which extends the first order case considered 
in \cite{ref7}(Theorem 3.2) to higher order cases. \cite{ref7} 
was later extended to the discrete time stochastic case by \cite{ref8}.
 In this paper, we aim to extend these results to deterministic 
higher order difference problems, using the ``squeezing'' argument. 

The application of higher order difference problems can be widely found in 
economics. In particular, they appear in the discussion concerning the 
overlapping generations models. 
 A satisfactory examination of 
the individuals' marriage and fertility decisions would necessitate the 
division of the representative agent's lifetime to multiple periods, instead 
of only two periods, young and old.
However, as argued in \cite{ref4}, the properties of a model with two-period-lived agents cannot be 
readily extended to $n$-period-lived agents.
 To consider the $n$-period-lived agents case, transversality 
conditions for higher order difference problems would be imperative. 

We first use the ``squeezing'' argument to derive both Euler's condition and 
the transversality condition for higher order difference problems, showing 
the argument needs two imperative assumptions. These two assumptions 
constitute the discrete time version of Assumption 1 and 2 in \cite{ref11}. 
We then provide a counterexample, in which the 
transversality condition is not satisfied without the two assumptions.
Because Assumption 1 and 2 are satisfied when a discounting 
factor is incorporated into the model, our transversality conditions also 
generalize the results obtained in the presence of discounting. For 
approaches on how to explicitly construct the optimal solutions to the undiscounted infinite horizon optimization
 problems, see \cite{ref2,ref3}.

\section{Derivation of the Transversality Conditions}

Suppose that the optimal path to (\ref{eq1}) exists and is given by ${\rm {\bf 
c}}^\ast\left( t \right) $, optimal in the sense of an overtaking criterion 
to be defined below. We perturb it with $N$th-order continuously differentiable curves 
${\rm {\bf q}}\left( t \right)$, 
\begin{equation}
\mathbf{c}(t)=\mathbf{c}^{\ast}(t)+\varepsilon \cdot \mathbf{q}(t).
\label{eq2}
\end{equation}
We define 
\begin{align}
 V(\varepsilon,T)=\displaystyle\inf_{T\leq T^{\prime}}\displaystyle\sum_{t=0}^{T^{\prime}}
 [&U(
 \mathbf{c}^{\ast}(t)+\varepsilon\cdot\mathbf{q}(t),
 \mathbf{c}^{\ast}(t+1)+\varepsilon\cdot\mathbf{q}(t+1),\nonumber\\
 &\qquad\cdots,
 \mathbf{c}^{\ast}(t+N-1)+\varepsilon\cdot\mathbf{q}(t+N-1),
 t)\nonumber\\
 &-U\left(\mathbf{c}^{\ast}(t),\mathbf{c}^{\ast}(t+1),\cdots,\mathbf{c}^{\ast}(t+N-1),t\right)].
 \label{eq3}
\end{align}
In this paper, \cite{ref1}'s notion of weak maximality is used as our optimality criterion.
We assume that there exists an optimal path that satisfy the weak 
maximality criterion, which is defined as: an attainable path $\left( {{\rm 
{\bf c}}^\ast \left( t \right)} \right)$ is optimal if no other attainable 
path overtakes it\footnote{ \cite{ref1} shows that such a path exists 
once two assumptions are satisfied. }:
\begin{align}
 \displaystyle\lim_{T\rightarrow+\infty}\displaystyle\inf_{T\leq T^{\prime}}\displaystyle\sum_{t=0}^{T^{\prime}}
 &[U(
 \mathbf{c}^{\ast}(t)+\varepsilon\cdot\mathbf{q}(t),
 \mathbf{c}^{\ast}(t+1)+\varepsilon\cdot\mathbf{q}(t+1),\nonumber\\
 &\qquad\cdots,
 \mathbf{c}^{\ast}(t+N-1)+\varepsilon\cdot\mathbf{q}(t+N-1),
 t)\nonumber\\
 &-U\left(\mathbf{c}^{\ast}(t),\mathbf{c}^{\ast}(t+1),\cdots,\mathbf{c}^{\ast}(t+N-1),t\right)]\leq 0.
 \label{eq4}
\end{align}

Let $V\left( \varepsilon \right)=\mathop {\lim }\limits_{T\to \infty } 
V\left( {\varepsilon ,T} \right)$. Differentiating it with respect to 
$\varepsilon $, we have
%
\begin{alignat}{2}
 &\displaystyle\lim_{\varepsilon\rightarrow {}^{+}0}\displaystyle\frac{V(\varepsilon)}{\varepsilon} \nonumber\\
=&\displaystyle\lim_{\varepsilon\rightarrow {}^{+}0}
\displaystyle\lim_{T\rightarrow\infty}\displaystyle\inf_{T\leq T^{\prime}}\displaystyle\sum_{t=0}^{T^{\prime}}
 &\displaystyle\frac{1}{\varepsilon}[U(
 \mathbf{c}^{\ast}(t)+\varepsilon\cdot\mathbf{q}(t),
 \mathbf{c}^{\ast}(t+1)+\varepsilon\cdot\mathbf{q}(t+1),\nonumber\\
 & &\qquad\cdots,
 \mathbf{c}^{\ast}(t+N-1)+\varepsilon\cdot\mathbf{q}(t+N-1),
 t)\nonumber\\
 & &-U\left(\mathbf{c}^{\ast}(t),\mathbf{c}^{\ast}(t+1),\cdots,\mathbf{c}^{\ast}(t+N-1),t\right)]. 
 \label{eq5}
\end{alignat}

%
%
Let $\displaystyle\lim_{\varepsilon\rightarrow{}^{+}0}
\displaystyle\frac{V(\varepsilon)}{\varepsilon}\equiv \Omega .$ 
Generally,
$\displaystyle\frac{d}{d\varepsilon}\displaystyle\lim_{T\rightarrow\infty}f(\varepsilon,T)
=\displaystyle\lim_{T\rightarrow\infty}\displaystyle\frac{d}{d\varepsilon}f(\varepsilon,T)$
 only when 
$\displaystyle\lim_{T\rightarrow\infty}\displaystyle\frac{d}{d\varepsilon}f(\varepsilon,T)$
 converges uniformly for $\varepsilon $ (\cite{ref9}).
We assume 


\textbf{Assumption 1.} Assume $\Omega$ converges uniformly for $\varepsilon $ when $T\rightarrow\infty$ .

Assume Assumption 1, we can then restate (\ref{eq5}) as
\begin{align}
 \Omega=&
 \displaystyle\lim_{T\rightarrow\infty}
 \displaystyle\lim_{\varepsilon\rightarrow{}^{+}0}
 \displaystyle\inf_{T\leq T^{\prime}}
 \displaystyle\frac{1}{\varepsilon}\nonumber\\
 &\times\displaystyle\sum_{t=0}^{T^{\prime}}[U\left(\mathbf{c}^{\ast}(t)+\varepsilon\cdot\mathbf{q}(t),\cdots,
 \mathbf{c}^{\ast}(t+N-1)+\varepsilon\cdot\mathbf{q}(t+N-1),t\right) \nonumber\\
 &\qquad\quad-U\left(\mathbf{c}^{\ast}(t),\cdots,
 \mathbf{c}^{\ast}(t+N-1),t\right)]. 
 \label{eq6}
\end{align}

We also assume 

\textbf{Assumption 2.}
 We assume for any 
$T>0,$
\begin{align*}
 \displaystyle\inf_{T\leq T^{\prime}}
 &\displaystyle\sum_{t=0}^{T^{\prime}}
 \displaystyle\frac{1}{\varepsilon}
 [U\left(\mathbf{c}^{\ast}(t)+\varepsilon\cdot\mathbf{q}(t),\cdots,
 \mathbf{c}^{\ast}(t+N-1)+\varepsilon\cdot\mathbf{q}(t+N-1),t\right) \nonumber\\
 &\qquad\quad-U\left(\mathbf{c}^{\ast}(t),\cdots,
 \mathbf{c}^{\ast}(t+N-1),t\right)] 
\end{align*}
%
%
converges uniformly for 
$\varepsilon .$


As in \cite{ref11}, a precise interpretation of Assumption 2 can be given as follows: Let 
\begin{align*}
 A(T,\varepsilon)=
 &\displaystyle\sum_{t=0}^{T^{\prime}}
 \displaystyle\frac{1}{\varepsilon}
 [U\left(\mathbf{c}^{\ast}(t)+\varepsilon\cdot\mathbf{q}(t),\cdots,
 \mathbf{c}^{\ast}(t+N-1)+\varepsilon\cdot\mathbf{q}(t+N-1),t\right) \nonumber\\
 &\qquad\quad-U\left(\mathbf{c}^{\ast}(t),\cdots,
 \mathbf{c}^{\ast}(t+N-1),t\right)]. 
\end{align*}
%
Then there exists a sequence $A\left(T_{n}^{\prime},\varepsilon\right)$ 
for each $\varepsilon>0$, so that 
$\displaystyle\lim_{n\rightarrow\infty}A\left(T_{n}^{\prime},\varepsilon\right)
=\displaystyle\inf_{T\leq T^{\prime}}A\left(T^{\prime},\varepsilon\right)$
%
, uniformly for $\varepsilon $, that is, the sequence 
is uniformly convergence for $\varepsilon $.

Assumptions 1 and 2 extend Assumption 3.1 in \cite{ref7}. When 
Assumptions 1 and 2 are satisfied, then $\mathop {\lim }\limits_{\varepsilon 
\to { }^+0} $ and $\mathop {\inf }\limits_{T\leqslant {T}'} $ can be 
interchanged, and equality (\ref{eq6}) can then restated as
\begin{align}
 \Omega=&
 \displaystyle\lim_{T\rightarrow\infty}
 \displaystyle\inf_{T\leq T^{\prime}}
 \displaystyle\lim_{\varepsilon\rightarrow{}^{+}0}
 \displaystyle\frac{1}{\varepsilon}\nonumber\\
 &\times\displaystyle\sum_{t=0}^{T^{\prime}}
 [U\left(\mathbf{c}^{\ast}(t)+\varepsilon\cdot\mathbf{q}(t),\cdots,
 \mathbf{c}^{\ast}(t+N-1)+\varepsilon\cdot\mathbf{q}(t+N-1),t\right) \nonumber\\
 &\qquad\quad-U\left(\mathbf{c}^{\ast}(t),\cdots,
 \mathbf{c}^{\ast}(t+N-1),t\right)]. 
 \label{eq7}
\end{align}
%
Because ${T}'$ is finite uniformly for $\varepsilon $, 
if 
\begin{align*}
 \displaystyle\sum_{t=0}^{T^{\prime}}
 &\displaystyle\frac{1}{\varepsilon}
 [U\left(\mathbf{c}^{\ast}(t)+\varepsilon\cdot\mathbf{q}(t),\cdots,
 \mathbf{c}^{\ast}(t+N-1)+\varepsilon\cdot\mathbf{q}(t+N-1),t\right) \nonumber\\
 &-U\left(\mathbf{c}^{\ast}(t),\cdots,
 \mathbf{c}^{\ast}(t+N-1)\right),t] 
\end{align*}
%
%
 exists, (\ref{eq7}) is then rewritten as
\begin{align}
 \Omega=&
 \displaystyle\lim_{T\rightarrow\infty}
 \displaystyle\inf_{T\leq T^{\prime}}
 \displaystyle\sum_{t=0}^{T^{\prime}}
 \displaystyle\frac{1}{\varepsilon}\nonumber\\
 &\times
 [U\left(\mathbf{c}^{\ast}(t)+\varepsilon\cdot\mathbf{q}(t),\cdots,
 \mathbf{c}^{\ast}(t+N-1)+\varepsilon\cdot\mathbf{q}(t+N-1),t\right) \nonumber\\
 &\qquad\quad-U\left(\mathbf{c}^{\ast}(t),\cdots,
 \mathbf{c}^{\ast}(t+N-1),t\right)]. 
 \label{eq8}
\end{align}
%
%
From the differentiability of $U$, we have

\begin{align*}
 & \displaystyle\lim_{\varepsilon\rightarrow{}^{+}0}
 \displaystyle\frac{1}{\varepsilon}
 [U\left(\mathbf{c}^{\ast}(t)+\varepsilon\cdot\mathbf{q}(t),\cdots,
 \mathbf{c}^{\ast}(t+N-1)+\varepsilon\cdot\mathbf{q}(t+N-1),t\right) \nonumber\\
 &\qquad\quad-U\left(\mathbf{c}^{\ast}(t),\cdots,
 \mathbf{c}^{\ast}(t+N-1),t\right)] \\ 
=&\displaystyle\sum_{i=1}^{n}
[\displaystyle\frac{\partial U\left(\mathbf{c}^{\ast}(t),\mathbf{c}^{\ast}(t+1),\cdots,\mathbf{c}^{\ast}(t+N-1),t\right)}
{\partial c_{i}(t)}q_{i}(t)\\
&\qquad+\displaystyle\frac{\partial U\left(\mathbf{c}^{\ast}(t),\mathbf{c}^{\ast}(t+1),\cdots,\mathbf{c}^{\ast}(t+N-1),t\right)}
{\partial c_{i}(t+1)}q_{i}(t+1)\\
&\qquad+\cdots
+\displaystyle\frac{\partial U\left(\mathbf{c}^{\ast}(t),\mathbf{c}^{\ast}(t+1),\cdots,\mathbf{c}^{\ast}(t+N-1),t\right)}
{\partial c_{i}(t+N-1)}q_{i}(t+N-1)].
%
\end{align*}
 Hence, 
\begin{align*}
\Omega=\displaystyle\lim_{T\rightarrow\infty}&\displaystyle\inf_{T\leq T^{\prime}}
\displaystyle\sum_{t=0}^{T^{\prime}}\displaystyle\sum_{i=1}^{n}
[\displaystyle\frac{\partial U\left(\mathbf{c}^{\ast}(t),\mathbf{c}^{\ast}(t+1),\cdots,\mathbf{c}^{\ast}(t+N-1),t\right)}
{\partial c_{i}(t)}q_{i}(t)
\nonumber\\
&\; +\displaystyle\frac{\partial U\left(\mathbf{c}^{\ast}(t),\mathbf{c}^{\ast}(t+1),\cdots,\mathbf{c}^{\ast}(t+N-1),t\right)}
{\partial c_{i}(t+1)}q_{i}(t+1)
+\cdots
\nonumber\\
&\; +\displaystyle\frac{\partial U\left(\mathbf{c}^{\ast}(t),\mathbf{c}^{\ast}(t+1),\cdots,\mathbf{c}^{\ast}(t+N-1),t\right)}
{\partial c_{i}(t+N-1)}q_{i}(t+N-1)].
\end{align*}
We derive
\begin{align}
&\displaystyle\sum_{t=0}^{T^{\prime}}\left(\displaystyle\sum_{i=1}^{n}\left(
\displaystyle\frac{\partial U(t)}{\partial c_{i}(t)}q_{i}(t)
+\cdots 
+\displaystyle\frac{\partial U(t)}{\partial c_{i}(t+N-1)}q_{i}(t+N-1)
\right)\right)\nonumber\\
=&\displaystyle\sum_{i=1}^{n}\{
\displaystyle\frac{\partial U(0)}{\partial c_{i}(0)}q_{i}(0)
+\displaystyle\frac{\partial \left(U(0)+U(1)\right)}{\partial c_{i}(1)}q_{i}(1)+\cdots\nonumber\\
&
+\displaystyle\frac{\partial\left(U(0)+\cdots+U(N-2)\right)}{\partial c_{i}(t)}q_{i}(N-2)\nonumber\\
&+\displaystyle\sum_{t=N-1}^{T^{\prime}}\displaystyle\frac{\partial\left(U(t-N+1)+\cdots+U(t)\right)}{\partial c_{i}(t)}q_{i}(t)
\nonumber\\
&+\left(\displaystyle\frac{\partial\left(U\left(T^{\prime}-N+2\right)+\cdots+U\left(T^{\prime}\right)\right)}
{\partial c_{i}(T^{\prime}+1)}\right)q_{i}\left(T^{\prime}+1\right)+\cdots\nonumber\\
&
+\displaystyle\frac{\partial U\left(T^{\prime}\right)}{\partial c_{i}\left(T^{\prime}+N-1\right)}q_{i}\left(T^{\prime}+N-1\right)\}.
\label{eq9}
\end{align}

Hence, Euler's condition is 
\begin{align}
&\frac{\partial U(0)}{\partial c_{i}(0)}=0,
\nonumber\\
&\frac{\partial \left(U(0)+U(1)\right)}{\partial c_{i}(1)}=0,
\nonumber\\
&\cdots\cdots,
\nonumber\\
&\frac{\partial \left(U(0)+\cdots+U(N-2)\right)}{\partial c_{i}(t)}=0,
\nonumber\\
%
%
%
&\frac{\partial \left(U(t-N+1)+\cdots+U(t)\right)}{\partial c_{i}(t)}=0,
\mbox{ for }
N-1\leq t\leq T^{\prime},
\label{eq10}
\end{align}

%
which extends the standard Euler's condition, 
and the transversality condition is
given by
\begin{align}
\displaystyle\lim_{T\rightarrow\infty}
\displaystyle\inf_{T\leq T^{\prime}}
\displaystyle\sum_{i=1}^{n}
[\displaystyle\frac
{\partial
\left(
U\left(T^{\prime}-N+2\right)+\cdots+U\left(T^{\prime}\right)
\right)
}
{\partial c_{i}\left(T^{\prime}+1\right)}
q_{i}\left(T^{\prime}+1\right)&
\nonumber\\
+\cdots+\displaystyle\frac{\partial\left(U\left(T^{\prime}\right)\right)}
{\partial c_{i}\left(T^{\prime}+N-1\right)}q_{i}\left(T^{\prime}+N-1\right)&]\leq 0.
\label{eq11}
\end{align}
Note that when $\varepsilon \to { }^-0$, the argument is the same:
\begin{align}
\displaystyle\lim_{T\rightarrow\infty}
\displaystyle\sup_{T\leq T^{\prime}}
\displaystyle\sum_{i=1}^{n}
[\displaystyle\frac
{\partial
\left(
U\left(T^{\prime}-N+2\right)+\cdots+U\left(T^{\prime}\right)
\right)
}
{\partial c_{i}\left(T^{\prime}+1\right)}
q_{i}\left(T^{\prime}+1\right)&
\nonumber\\
+\cdots+\displaystyle\frac{\partial\left(U\left(T^{\prime}\right)\right)}
{\partial c_{i}\left(T^{\prime}+N-1\right)}q_{i}\left(T^{\prime}+N-1\right)&]\geq 0.
\tag{$11^{\prime}$}
\label{eq12}
\end{align}
%
%
Next, 
we consider the linkage between our result and that in \cite{ref7}.
We fix $0<\bar{\alpha}<1$ 
and 
$\alpha:\mathbb{R}^{+}\rightarrow\mathbb{R}^{+}$,
$C^{\infty}$, 
$\alpha(0)=0$,
$\cdots$,
$\alpha^{(n-1)}(0)=0$,
$\alpha(t)=\bar{\alpha}$,
$t\geq 1$.
 We let $\varepsilon \to { }^+0$. 
Let $q(t)=\alpha{c}^{\ast}(t)$, 
then (\ref{eq11}) is modified to 
%
\begin{align}
&\displaystyle\lim_{T\rightarrow\infty}
\displaystyle\inf_{T\leq T^{\prime}}
\displaystyle\sum_{i=1}^{n}
[\displaystyle\frac
{\partial
\left(
U\left(T^{\prime}-N+2\right)+\cdots+U\left(T^{\prime}\right)
\right)
}
{\partial c_{i}\left(T^{\prime}+1\right)}
\alpha c_{i}^{\ast}\left(T^{\prime}+1\right)
\nonumber\\
&\qquad\qquad\qquad\quad+\cdots+\displaystyle\frac{\partial\left(U\left(T^{\prime}\right)\right)}
{\partial c_{i}\left(T^{\prime}+N-1\right)}\alpha c_{i}^{\ast}\left(T^{\prime}+N-1\right)]\nonumber\\
=&\bar{\alpha}\displaystyle\lim_{T\rightarrow\infty}
\displaystyle\inf_{T\leq T^{\prime}}
\displaystyle\sum_{i=1}^{n}
[\displaystyle\frac
{\partial
\left(
U\left(T^{\prime}-N+2\right)+\cdots+U\left(T^{\prime}\right)
\right)
}
{\partial c_{i}\left(T^{\prime}+1\right)}
c_{i}^{\ast}\left(T^{\prime}+1\right)
\nonumber\\
&\qquad\qquad\qquad\quad+\cdots+\displaystyle\frac{\partial\left(U\left(T^{\prime}\right)\right)}
{\partial c_{i}\left(T^{\prime}+N-1\right)}
c_{i}^{\ast}\left(T^{\prime}+N-1\right)]\nonumber\\
\leq &0.
\label{eq13}
\end{align}
%
Because $\bar {\alpha }>0$, we then have
\begin{align}
\displaystyle\lim_{T\rightarrow\infty}
\displaystyle\inf_{T\leq T^{\prime}}
\displaystyle\sum_{i=1}^{n}
[\displaystyle\frac
{\partial
\left(
U\left(T^{\prime}-N+2\right)+\cdots+U\left(T^{\prime}\right)
\right)
}
{\partial c_{i}\left(T^{\prime}+1\right)}
c_{i}^{\ast}\left(T^{\prime}+1\right)&
\nonumber\\
\qquad\qquad\qquad\quad+\cdots+\displaystyle\frac{\partial\left(U\left(T^{\prime}\right)\right)}
{\partial c_{i}\left(T^{\prime}+N-1\right)}
c_{i}^{\ast}\left(T^{\prime}+N-1\right)]&\leq 0.
\label{eq14}
\end{align}
which is an extension of \cite{ref7}'s transversality condition. 

\section{A Counterexample}

We proceed to show that Assumption 1 and 2 
are imperative in the sense that (\ref{eq11}) becomes invalid if one of them is 
violated. We consider the following simple counterexample:
\begin{equation}
 U\left(c(t),c(t+1),c(t+2),t\right)
=\left(c(t)-\alpha\right)^{2}+\beta c(t+1)+\gamma c(t+2),
\label{eq15}
\end{equation}
where $\alpha>0$, 
$\beta>0$, 
$\gamma>0$, 
and the initial values 
$c(0)=c_{0}$,
$c(1)=c_{1}$ 
are given.
From (\ref{eq10}), we see that Euler's condition is given by
\begin{align}
&\frac{\partial U(0)}{\partial c(0)}=0,
\nonumber\\
&\frac{\partial \left(U(0)+U(1)\right)}{\partial c(1)}=0,
\nonumber\\
&\frac{\partial \left(U(t-2)+U(t-1)+U(t)\right)}{\partial c(t)}=0,
2\leq t\leq {T}^\prime,
\label{eq16}
\end{align}
which implies
\begin{align}
 t=0,\quad&2\left(c(0)-\alpha\right)=0, \tag{$15^{\prime}$}\\
 t=1,\quad&2\left(c(1)-\alpha\right)+\beta=0, \tag{$15^{''}$}\\
 t=2,\quad&2\left(c(2)-\alpha\right)+\beta+\gamma=0, \tag{$15^{'''}$}\\
 t=3,\quad&2\left(c(3)-\alpha\right)+\beta+\gamma=0, \tag{$15^{''''}$} \\
 \cdots\cdots,& \nonumber\\
 t=T^{\prime},\quad&2\left(c\left(T^{\prime}\right)-\alpha\right)+\beta+\gamma=0.\tag{$15^{'''''}$} 
\end{align}
Thus, we have 
$c(2)
=c(3)
=\cdots
=c\left(T^{\prime}\right)
=\alpha-\displaystyle\frac{\beta+\gamma}{2}$. 

Choosing a $p$ so that $p(0)=0$ and $p(t)>0$, 
there exists $T_{0}>0$, 
$p(t)$ is a constant $p_{\infty}>0$ 
when $t\geq T_{0}$. 

From (\ref{eq15}), 
we see that
\begin{align}
 &\displaystyle\frac{\partial\left(U\left(T^{\prime}-1\right)+U\left(T^{\prime}\right)\right)}{\partial c(t+1)}q\left(T^{\prime}+1\right)
 =(\gamma+\beta)q\left(T^{\prime}+1\right)\leq 0,
 \label{eq20} \\
 &\displaystyle\frac{\partial\left(U\left(T^{\prime}\right)\right)}{\partial c(t+2)}q\left(T^{\prime}+2\right)
 =\gamma q\left(T^{\prime}+2\right)\leq 0.
 \label{eq21} 
\end{align}
Hence, 
we have arrived at a contradiction to (\ref{eq11}). 


Next, we show that Assumption 1 is violated, which causes this contradiction. 
We consider 
%
%
$U\left(
c^{\ast}(t)+\varepsilon p(t),
c^{\ast}(t+1)+\varepsilon p(t+1),
c^{\ast}(t+2)+\varepsilon p(t+2)\right)$
$-U\left(c^{\ast}(t),c^{\ast}(t+1),c^{\ast}(t+2)\right)$.
%
Substituting $c^{\ast}(t)=\alpha-\displaystyle\frac{\beta+\gamma}{2}$ into it, 
we have
\begin{align}
 &U\left(
c^{\ast}(t)+\varepsilon p(t),
c^{\ast}(t+1)+\varepsilon p(t+1),
c^{\ast}(t+2)+\varepsilon p(t+2)\right) \nonumber\\
 &-U\left(c^{\ast}(t),c^{\ast}(t+1),c^{\ast}(t+2)\right) \nonumber \\
=&\left(c^{\ast}(t)+\varepsilon p(t)-\alpha\right)^{2}
 +\beta\left(c^{\ast}(t+1)+\varepsilon p(t+1)\right)
 +\gamma\left(c^{\ast}(t+2)+\varepsilon p(t+2)\right) \nonumber \\
 &-\left(\left(c^{\ast}(t)-\alpha\right)^{2}+\beta c^{\ast}(t+1)+\gamma c^{\ast}(t+2)\right) \nonumber \\
=&\left(\varepsilon p(t)-\left(\displaystyle\frac{\beta+\gamma}{2}\right)\right)^{2}
 +\varepsilon\left(\beta p(t+1)+\gamma p(t+2)\right)
 -\left(\displaystyle\frac{\beta+\gamma}{2}\right)^{2}.
\label{eq22}
\end{align}
Hence, 
\begin{align}
 &\displaystyle\inf_{T\leq T^{\prime}}
 \displaystyle\sum_{t=0}^{T^{\prime}}
 \displaystyle\frac{\left(\varepsilon p(t)
-\left(\displaystyle\frac{\beta+\gamma}{2}\right)\right)^{2}
+\varepsilon \beta p(t+1)+\varepsilon \gamma p(t+2)
-\left(\displaystyle\frac{\beta+\gamma}{2}\right)^{2}}{\varepsilon}
\nonumber \\ 
=&\displaystyle\inf_{T\leq T^{\prime}}
 \displaystyle\sum_{t=0}^{T^{\prime}}
 \left(
 \varepsilon p(t)^{2}-(\beta+\gamma)p(t)
 +\beta p(t+1)
 +\gamma p(t+2)
 \right) 
\nonumber\\ 
 =&\displaystyle\inf_{T\leq T^{\prime}}
 \left(\varepsilon 
\displaystyle\sum_{t=0}^{T^{\prime}}p(t)^{2}
+\beta p\left(T^{\prime}+1\right)
+\gamma p\left(T^{\prime}+2\right)
\right) \nonumber\\ 
 =&\displaystyle\inf_{T\leq T^{\prime}}\left(
 \varepsilon \displaystyle\sum_{t=0}^{T^{\prime}}\left(p(t)^{2}\right)
\right)
+\beta p_{\infty}
+\gamma p_{\infty} 
\nonumber\\ 
 =&\varepsilon \displaystyle\sum_{t=0}^{T}\left(p(t)^{2}\right) 
+\beta p_{\infty} 
+\gamma p_{\infty}. 
\label{eq23}
\end{align}

$\Omega$ is the limit of (\ref{eq22}) 
when $T\rightarrow\infty$, 
$\varepsilon\rightarrow 0$. 
However, \\
because 
$
\displaystyle\lim_{\varepsilon\rightarrow 0}
\displaystyle\lim_{T\rightarrow\infty}
\left(\varepsilon\displaystyle\sum_{t=0}^{T}\left(p(t)^{2}\right)+\beta p_{\infty}+\gamma p_{\infty}\right)=\infty,$\\
whereas 
$
\displaystyle\lim_{T\rightarrow\infty}
\displaystyle\lim_{\varepsilon\rightarrow 0}
\left(\varepsilon\displaystyle\sum_{t=0}^{T}\left(p(t)^{2}\right)+\beta p_{\infty}+\gamma p_{\infty}\right)
=\beta p_{\infty}+\gamma p_{\infty},$
we see that $\Omega $ does not converge uniformly 
for $\varepsilon$ when $T\rightarrow\infty$. 
Hence, 
Assumption 1 is violated 
and (\ref{eq11}) is also not satisfied.

\section{Conclusion}

This paper gives the two assumptions that would be imperative when examining 
infinite horizon discrete time optimization problems in which the objective functions are 
unbounded. Our results generalizes the results of \cite{ref5,ref6,ref7,ref10,ref12}($N=1)$ to higher order difference 
problems. Specifically, when $N=1$, our transversality condition is the 
discrete time version of that of \cite{ref7}(Theorem 3.2).
 Moreover, 
our Assumption 1 and 2 also constitute the discrete time version of the Assumptions in \cite{ref11}.
As in \cite{ref11},
 our Assumption 1 and 2 obviously hold when a discounting 
factor is incorporated into the model. In this sense, this paper also 
extends the transversality conditions examined in the presence of 
discounting. 

\medskip


\label{lastpage}
\end{document}